\documentclass[12pt]{iopart}

\usepackage{mathrsfs,amssymb,amsthm}

\newtheorem{theorem}{Theorem}[section]
\newtheorem{lemma}[theorem]{Lemma}

\newcommand{\ud}{\mathrm{d}}
\newcommand{\R}{\mathbb{R}}
\newcommand{\cA}{\mathcal{A}}

\begin{document}

\title[A simple proof of uniqueness of trajectories]{A simple proof of uniqueness of the particle trajectories for solutions of the Navier-Stokes equations}

\author{M Dashti, J C Robinson}

\address{Mathematics Institute, University of Warwick, Coventry CV4 7AL, UK}
\ead{M.Dashti@warwick.ac.uk, J.C.Robinson@warwick.ac.uk}

\begin{abstract}
We give a simple proof of the uniqueness of fluid particle trajectories corresponding to: 1) the solution of the two-dimensional Navier Stokes equations with an initial condition that is only square integrable, and 2) the local strong solution of the three-dimensional equations with an $H^{1/2}$-regular initial condition i.e.\ with the minimal Sobolev regularity known to guarantee uniqueness. This result was proved by Chemin \& Lerner (J Diff Eq 121 (1995) 314-328) using the Littlewood-Paley theory for the flow in the whole space $\R^d$, $d\ge 2$. 
We first show that the solutions of the differential equation $\dot{X}=u(X,t)$ are unique if $u\in L^p(0,T;H^{(d/2)-1})$ for some $p>1$ and $\sqrt{t}\,u\in L^2(0,T;H^{(d/2)+1})$. We then prove, using standard energy methods, that the solution of the Navier-Stokes equations with initial condition in $H^{(d/2)-1}$ satisfies these conditions. This proof is also valid for the more physically relevant case of bounded domains.
\end{abstract}

\submitto{\NL}
\maketitle

\section{Introduction}\label{sec:intro}

We consider the Navier-Stokes equations
\begin{equation}\label{eq:nse}
\frac{\partial u}{\partial t}-\nu\Delta u+(u\cdot \nabla)u+\nabla p=f,\qquad u(x,0)=u_0,
\quad u|_{\partial\Omega}=0
\end{equation}
in which $x\in\Omega\subset \R^d$ with $d=2,3$, $\Omega$ is an open set with a sufficiently smooth boundary,
$u(x,t)$ is the velocity vector field,
$p(x,t)$ the pressure scalar function,
$f(x,t)$ the body force,
and $\nu$ is the kinematic viscosity which is considered constant.

The minimal Sobolev regularity for the initial condition that is known to give rise to a unique solution 
$u\in L^\infty(0,T;H^{(d/2)-1}(\Omega))\cap L^2(0,T;H^{d/2}(\Omega))$ is $u_0\in L^2(\Omega)$ 
for the two-dimensional domain (Leray 1933 for the whole plane, Lions \& Prodi 1959 and Ladyzhenskaya 1958 for bounded domains) 
and $u_0\in H^{1/2}(\Omega)$ in the case of a three-dimensional domain (Fujita \& Kato 1964), where $H^s$ with real $s>0$ is the standard Sobolev space of order $s$ (we recall the characterization of fractional Sobolev spaces in Section \ref{sec:3d}).
For the two-dimensional Navier-Stokes equations the above unique solution is global in time, while in
the three-dimensional case a unique solution exists on $[0,T_1)$ and the best available bound on
the $H^{1/2}$-norm of $u(t)$ tends to infinity as $t\to T_1$. 
We consider $T<T_1$, in the three-dimensional case.
In this paper we denote by $(\cdot,\cdot)$ and $\|\cdot\|$ the inner product and norm in $L^2(\Omega)$, respectively.
The norm in $H^r(\Omega)$ is denoted by $\|\cdot\|_r$ and the norm in any other normed space $E$, by $\|\cdot\|_E$. 

Corresponding to the solution $u$ defined above, as long as it exists, the fluid particle trajectories are the solutions of
\begin{equation}\label{eq:ode}
\frac{\ud X}{\ud t}=u(X,t),\qquad X(0)=a\in\Omega.
\end{equation}
At least one solution to the above system exists. 
This can be shown (following Foias, Guillop\'e \& Temam 1985) by considering $u_n$ to be the Galerkin approximations of $u$ and defining $X_n$ to be the solution of 
\begin{equation*}
\frac{\ud X_n}{\ud t}=u_n(X_n,t),\quad\mbox{with }X_n(0)=a
\end{equation*}
and then showing the uniform convergence of $X_n$ to $X$ in $[0,T]$ and strong convergence of $u_n$ to $u$ in $L^1(0,T;L^\infty)$. We adapt and explain this argument in the proof of the existence of solutions in Theorem \ref{thm:g-uniqueness}.

The uniqueness of the solutions of (\ref{eq:ode}) in the whole space is shown by Chemin \& Lerner (1995) using the Littlewood-Paley theory. In this paper we present an alternative simpler proof which is valid in a general bounded domain as well.

%
%

\subsection{Chemin and Lerner's proof of uniqueness}

The uniqueness of the solutions of (\ref{eq:ode}) is shown 
by Chemin and Lerner in their 1995 paper. They use the Littlewood-Paley theory to prove 
enough regularity for the solution of the Navier-Stokes equations in order to be able to apply a 
generalisation of the Osgood criterion to (\ref{eq:ode}). 
They denote by $\mathscr{H}_{1,T}^{(d/2)+1}$ the space of functions $u:[0,T]\times\R^d\to \R^d$ 
satisfying
$$
\left( \sum_{q\in \mathbb{N}}2^{q(2+d)} \left( \int_0^T \|\Delta_q u(t)\|_{L^2} \right)^2 \right)^{1/2}<\infty
$$
where
$$
\Delta_q u=2^{qd}\int_{\R^d} h(2^q y)u(x-y)\,\ud y,
$$
and $h$ is the inverse fourier transform of some $\phi$ that is an appropriate bump function supported 
on the annulus $\{3/4\le|\xi|\le8/3\}$.
They prove that the solution of the two-dimensional Navier-Stokes 
equations is an element of $\mathscr{H}_{1,T}^{(d/2)+1}$  and then show that 
$\mathscr{H}_{1,T}^{(d/2)+1}\subset L^1_\mathrm{loc}(0,T;C_{\omega_{\epsilon}}(\R^d;\R^d))$ 
where $\omega_\epsilon(r)=r(1-\log r)^{\epsilon+1/2}$ and $C_{\omega}(\R^d;\R^d)$ is a Banach space with the norm
$$
\|u\|_{\omega}=\|u\|_{L^\infty(\R^d)}+\sup_{x,y\in \R^d\times \R^d,\, x\neq y}\frac{|u(x)-u(y)|}{\omega(|x-y|)}.
$$
They conclude the uniqueness of the flow by proving a generalisation of the Osgood criterion
which states that if $F\in L_{\mathrm{loc}}^1(0,T,C_{\omega})$ with $\omega$ satisfying
$$
\int_0^1 \frac{\ud r}{\omega(r)}=+\infty,
$$ 
then the equation 
$$
x(t)=x_0+\int_0^t F(s,x(s))\,\ud s
$$
has a unique solution over $[0,t_1]$ for some $t_1<T$.\\

%
%

\subsection{The summary of our proof}

We present an alternative proof of the same uniqueness result in the case of bounded two- and 
three-dimensional domains, which is also valid for the whole space $\R^d$ and periodic domains.

The proof is in fact elementary. 
We denote by $\eta(t)$ the Euclidean norm of the difference of 
two solutions of (\ref{eq:ode}) at time $t$, write down the differential equation satisfied by $\eta(t)$, and derive an upper bound for this difference
in terms of the vector field $u$ and the value of $\eta$ at some previous time $s>0$.
Letting $s\to 0$ however, is not straightforward. We outline the difficulty here and address it in the following sections.

Assuming that both $X(t)$ and $Y(t)$ satisfy (\ref{eq:ode}), we have
$$
\frac{\ud}{\ud t}(X-Y)=u(X,t)-u(Y,t),\qquad \mbox{with } X(0)-Y(0)=0.
$$
Using a result of Zuazua (2002), we know that a vector field $u$ over a $d$-dimensional domain $\Omega$, satisfies
$$
|u(X)-u(Y)|\le c\|u\|_{{1+d/2}}\,|X-Y|\,(-\log{|X-Y|})^{1/2}
$$
for any $X,Y\in \Omega$. This bound is obtained easily by considering the extension of $u$ to $\R^d$, $E[u]$, and writing $E[u]$ as the inverse of its Fourier 
transform.
Once we have the above inequality we can write, for $\eta(t)=|X(t)-Y(t)|$,
$$
\frac{\ud\eta}{\ud t} \le c\|u\|_{{1+d/2}}\,\eta\,(-\log\eta)^{1/2} ,
$$
implying that
\begin{equation}\label{eq:eta}
\eta(t) \,\le\,\exp \left(-\left( \big(\log({1}/{\eta(s)})\,\big)^{1/2}
           -c\int_s^t\|u\|_{{1+d/2}}\,\ud\tau \right)^2\right)
\end{equation}
for $0<s<t$. The uniqueness of the solutions of (\ref{eq:ode}) follows by showing that the right-hand side of the above inequality converges to zero as $s\to 0$. If $\int_0^t\|u\|_{1+d/2}\,\ud\tau<\infty$, then one could simply let $s\to 0$ to obtain the result. The problem is that it is not known whether $u\in L^1(0,T;H^{1+d/2})$
(in fact even for the heat equation this is not known to be true). To circumvent this,
we bound $\int_s^t\|u\|_{{1+d/2}}\,\ud\tau$ and $\eta(s)$ from above and show that
$$
\lim_{s\to 0}\left( \big(\log({1}/{\eta(s)})\,\big)^{1/2}
           -c\int_s^t\|u\|_{{1+d/2}}\,\ud\tau\right) =\infty
$$
for small enough $t$. We note that $u$ is smooth for any $t>0$ and therefore showing the uniqueness for an arbitrary small interval containing $t=0$ is enough.

In Section \ref{sec:g-uniqueness} we give sufficient conditions on the vector function $u$
that result in appropriate bounds on the logarithmic and integral terms in the 
right-hand side of (\ref{eq:eta}), as discussed above, to ensure the uniqueness of the solution of the ordinary differential equation (\ref{eq:ode}). We require that $u$ satisfies
\begin{equation}\label{eq:mcond0}
\fl\qquad u\in L^p(0,T;H^{(d/2)-1}(\Omega))\quad\mbox{with } p>1,
\quad\mbox{ and }\quad \sqrt{t}\,u\in L^2(0,T;H^{1+d/2}(\Omega)).
\end{equation}

We then in Section \ref{sec:estimates}, show that the solution of the Navier-Stokes equations (\ref{eq:nse}) with $u_0\in H^{(d/2)-1}(\Omega)$ and $f\in L^2(0,T; H^{(d/2)-1}(\Omega))$ satisfies (\ref{eq:mcond0}). 
This is straightforward in the two-dimensional case. 
For bounded three-dimensional domains, we need to use 
the fractional powers of $\cA=-\Delta$ with $D(\cA)=\{u\in H^2(\Omega): u|_{\partial\Omega}=0\}$. 
To deal with these, we  need the equivalence of
$\|\cA^{r/2}u\|$ and $\|u\|_r$ for any non-negative real number $r$
and any $u\in D(\cA^{r/2})$. 
For $0\le r\le 2$, a concrete characterisation of $D(\cA^{r/2})$ is known 
(Fujiwara 1967) which gives the equivalence of
$\|\cA^{r/2}u\|$ and $\|u\|_r$ for $0\le r\le 2$. 
But, to our knowledge, such a characterisation does not exist
for $r>2$. 
The equivalence of the norms in $D(\cA^{r/2})$ and $H^r$ for any $r\in\R_{\ge 0}$ however, 
can be concluded almost immediately from an interpolation theorem of Lions \& Magenes (1972). 
We show this in Section \ref{sec:3d} and use it to prove the validity of (\ref{eq:mcond0}) for the solutions of the three-dimensional Navier-Stokes equations.

%
%
\section{A general uniqueness result}\label{sec:g-uniqueness}

In this section we prove a uniqueness result for general ordinary differential
equations. It will be then shown in the next section, that the uniqueness condition of this 
theorem is satisfied by the solution of the Navier-Stokes equations.

\begin{theorem}\label{thm:g-uniqueness}

Let $\Omega$ be the whole space $\R^d$, $d\ge 2$, a periodic $d$-dimensional domain or an open bounded subset of $\R^d$ 
with a sufficiently smooth boundary. Consider 
\begin{equation}\label{eq:mcond}
\fl\qquad u\in L^p(0,T;H^{(d/2)-1}(\Omega))\quad\mbox{with } p>1,
\quad\mbox{ and }\quad \sqrt{t}\,u\in L^2(0,T;H^{(d/2)+1}(\Omega)),
\end{equation}
with $u=0$ on $\partial\Omega$ when $\Omega$ is a bounded domain. Then the differential equation
\begin{equation}\label{eq:ode2}
\frac{\ud X(t)}{\ud t}=u(X(t),t),\qquad X(0)=a\in\Omega
\end{equation}
has a unique solution over $[0,T]$. 
\end{theorem}


We note that in the case of a bounded domain $\Omega$, assuming $\partial\Omega$ to satisfy the uniform $C^1$-regularity condition (Adams 1975) is sufficient.  This ensures that the trace operator, mapping $u$ to $u|_{\partial\Omega}$, makes sense (Adams 1975, Theorem 7.53), and also the extension operator used to derive the Sobolev embedding results on $\Omega$ from those on $\R^d$ is well-defined (Adams 1975, Theorem 4.32).

\begin{proof}

To prove the existence, we first show the integrability of $\|u\|_{L^\infty(\Omega)}$ over 
$(0,T)$. By Agmon's inequality we have
$$
\|u\|_{L^\infty(\Omega)}\,\le\, c\,\|u\|_{(d/2)-1}^{1/2}\,\|u\|_{(d/2)+1}^{1/2}
$$
and therefore we can write
\begin{eqnarray}\label{eq:L1-infty}
&\fl\qquad\int_0^t\|u\|_{L^\infty(\Omega)}\,\ud\tau
\le c\,\int_0^t \tau^{-1/4}\,\|u\|_{(d/2)-1}^{1/2}\,\tau^{1/4}\,\|u\|_{(d/2)+1}^{1/2}\,\ud \tau
                    \nonumber\\
&\fl\qquad\qquad\le c\left( \int_0^t \tau^{-p/(3p-2)}\,\ud \tau \right)^\frac{3p-2}{4p}
         \left( \int_0^t \|u\|_{(d/2)-1}^p\,\ud \tau \right)^\frac{1}{2p}
         \left( \int_0^t t\|u\|_{(d/2)+1}^2\,\ud \tau \right)^\frac{1}{4}\nonumber\\
&\fl\qquad\qquad\le\, c\,t^{\frac{p-1}{2p}}.
\end{eqnarray}
For a bounded or periodic domain $\Omega$, we let $u_n$ be the $n$-dimensional Galerkin approximation of $u$ (the image of $u$ under the projection in $L^2$ onto the space spanned by the first $n$ eigenfunctions of $\cA=-\Delta$ with Dirichlet boundary conditions when $\Omega$ is bounded and periodic boundary conditions for a periodic domain). For the whole space $\R^d$, we consider $u_n$ to be a sequence of mollified versions of $u$ with $1/n$ the parameter of mollification. Following Foias, Guillop\'e \& Temam (1985) we consider $X_n$ to be the solution of
$$
\frac{\ud X_n}{\ud t}=u_n(X_n,t),\quad\mbox{with }X_n(0)=a.
$$
Now, $u_n$ is continuous in $X_n$ and also has the same regularity properties of $u$ and therefore is integrable with respect to $t$.
Hence a solution of the above differential system exists over $[0,T]$ unless, in the case of $\Omega$ a bounded domain, the particle leaves $\Omega$ at some time less than $T$.
This is not possible  since $u_n$ is zero on $\partial\Omega$ and hence (\ref{eq:ode}) is
solvable over $[0,T]$.
By (\ref{eq:L1-infty}), we also have
$$
|X_n(t)|\le \int_0^t |u_n(X_n,t)|\,\ud \tau<c
$$
and 
$$
|X_n(t)-X_n(s)|\le \int_s^t \|u(\tau)\|_{L^\infty(\Omega)}\,\ud\tau \le c\,|t-s|^{\frac{p-1}{2p}}
$$
implying that $\{X_n\}$ is equicontinuous. Therefore, by the Arzel\`a-Ascoli Theorem, $X_n$ has a subsequence (which we label by $n$ again) that converges uniformly
to some $X(t)$. 
To prove that $X(t)$ is the solution of (\ref{eq:ode2}) we need to show that $u_n(X_n,t)\to u(X,t)$ in $L^1(0,T)$. We can write
\begin{eqnarray}\label{eq:galerkin}
\int_0^T |u_n(X_n,t)- u(X,t)|\,\ud t
&&\le \int_0^T |u_n(X_n,t)- u(X_n,t)|\,\ud t\nonumber\\
&&\quad+\int_0^T |u(X_n,t)- u(X,t)|\,\ud t.
\end{eqnarray}
The first term in the right-hand side converges to zero since we know that
$u_n\to u$ in $L^1(0,T;L^\infty)$.
For the second term in the right-hand side of (\ref{eq:galerkin}) we note that for almost every strictly positive $t$,
$u(t)\in H^{(d/2)+1}$  and therefore is continuous in 
$x$, implying that $u(X_n,t)\to u(X,t)$ as $n\to \infty$ for almost every $t\in (0,T)$. Since 
$u(x,t)$ is integrable over $(0,T)$, the result follows by the Lebesgue Dominated Convergence Theorem.\\

We now prove uniqueness.
In Section \ref{sec:intro}, we showed that $\eta(t)=|X(t)-Y(t)|$ satisfies
\begin{equation}\label{eq:eta2}
\eta(t) \,\le\,\exp \left(-\left( \big(\log({1}/{\eta(s)})\,\big)^{1/2}
           -c\int_s^t\|u\|_{{1+d/2}}\,\ud\tau \right)^2\right)
\end{equation}
for $0<s<t$. We need to obtain appropriate lower and upper bounds for the logarithmic and integral term respectively in the right-hand side of the above inequality as $s\to 0$. 
To bound $\log(1/\eta(s))$, we write
\begin{equation*}
\frac{\ud \eta}{\ud t}\le |u(X,t)-u(Y,t)|
\le 2\|u\|_{L^\infty(\Omega)}.
\end{equation*}
By (\ref{eq:L1-infty}) we have
$$
\eta(s)\le c\int_0^s \|u\|_{L^\infty(\Omega)}\,\ud\tau
\le c\,s^{\frac{p-1}{2p}}
$$
which implies that
$$
\log \frac{1}{\eta(s)}\,\ge\, c\,\log\frac{1}{s}.
$$
for $s$ sufficiently small.

For the integral term in the right-hand side of (\ref{eq:eta2}) we write
\begin{eqnarray*}
\int_s^t \|u\|_{1+d/2}\,\ud\tau
&\le\int_s^t \tau^{-1/2}\,\tau^{1/2}\,\|u\|_{1+d/2}\,\ud\tau\\
&\le \left(\int_s^t \tau^{-1}\,\ud\tau \right)^{1/2}\,
               \left( \int_s^t \tau\|u\|_{1+d/2}^2\,\ud\tau \right)^{1/2}\\
&\le c\,(\log t-\log s)^{1/2}\,\left( \int_s^t \tau\|u\|_{1+d/2}^2\,\ud\tau \right)^{1/2}
\end{eqnarray*}
and therefore for $s$ small enough
$$
\int_s^t \|u\|_{1+d/2}\,\ud\tau \,\le\, c\,K_d(t) \left(\log t+ \log\frac{1}{s} \right)^{1/2}
$$
with
$$
K_d^2(t)=\int_0^t \tau\|u\|_{1+d/2}^2\,\ud\tau.
$$

Having the above bounds, we go back to (\ref{eq:eta2}), fix $t$ and let $s\to 0$ in the right-hand side to write
\begin{eqnarray*}
\eta(t)\,
\le\,\exp \left( -\lim_{s\to 0}\log(1/s)\,(1-c\,K_d(t))^2\right).
\end{eqnarray*}
Since $K_d^2(t)$ is the integral of an integrable function, 
it is absolutely continuous (Priestley 1997) and therefore we can choose $T^*$ small enough so that
$$
1-c\,K_d(T^*) > 0.
$$
Hence $\eta(t)\to 0$ for any $t\le T^*$ as $s\to 0$. This gives the result for $t\in[0,T^*]$. The uniqueness over $[T^*, T]$ then follows easily from (\ref{eq:eta2}) since $u\in L^1(T^*,T;H^{(d/2)+1})$.
\end{proof}
\vskip.3cm

We note that the assumption $\sqrt{t}\,u\in L^2(0,T; H^{(d/2)+1})$ implies that $u\in L^r(0,T; H^{(d/2)+1})$ for any $r<1$.\vskip.3cm


Using an argument similar to the proof of the existence part of the above theorem, one can show that if a sequence $\sqrt{t}\, u_n$ is uniformly bounded in $L^2(0,T;H^{(d/2)+1}(\Omega))$ and $u_n\to u$ strongly in $L^p(0,T;H^{(d/2)-1}(\Omega))$ with $p>1$, then the solution of $\dot{X}_n=u_n(X_n,t)$ converges uniformly to the solution of $\dot{X}=u(X,t)$ (Dashti \& Robinson 2008). Using this in the case of the Navier-Stokes equations, one can show that the map $u_0 \mapsto X(t)$ is continuous from $H^{(d/2)-1}(\Omega)$ into $C^0([0,T);\R^d)$.\vskip.3cm

%
%
\section{The estimates for the solutions of the Navier-Stokes equations}\label{sec:estimates}

In this section we show that $u$, a unique solution of the Navier-Stokes equations (\ref{eq:nse}) with $u_0\in H^{(d/2)-1}(\Omega)$, 
satisfies the conditions of Theorem \ref{thm:g-uniqueness} and thereby prove the following theorem:

\begin{theorem}\label{thm:main}
Let $\Omega$ be an open bounded subset of $\R^d$, $d= 2,3$, with a sufficiently smooth boundary.
Consider $u\in L^\infty(0,T;H^{(d/2)-1}(\Omega))\cap L^2(0,T;H^{d/2}(\Omega))$ to be a unique solution of 
the Navier-Stokes equations (\ref{eq:nse}) with 
$u_0\in H^{(d/2)-1}(\Omega)$ and 
$f\in L^2(0,T;H^{(d/2)-1}(\Omega))$.
Then the ordinary differential system
$$
\frac{\ud X(t)}{\ud t}=u(X(t),t),\quad\mbox{ with}\quad X(0)=a,
$$
has a unique solution.
\end{theorem}

We note that for the above result to be true, assuming $\partial\Omega$ to satisfy uniform $C^1$-regularity condition is sufficient. For such a domain both the boundary trace embedding theorem (Adams 1975, Theorem 5.22) that we need for obtaining the estimates of Lemma \ref{lem:2d} and \ref{lem:3d}, and the result of Theorem \ref{thm:g-uniqueness} hold.


\subsection{The two-dimensional case}

The following result is due to Temam (1977). 

\begin{lemma}\label{lem:2d}
Let $\Omega$ be an open bounded subset of $\R^2$ with a sufficiently smooth boundary.
Consider $u\in L^\infty(0,T;L^2(\Omega))\cap L^2(0,T;H^{1}(\Omega))$ to be a unique solution of 
the two-dimensional Navier-Stokes equations (\ref{eq:nse}) with
$u_0\in L^2(\Omega)$ and $f\in L^2(0,T;L^2(\Omega))$. 
Then 
\begin{eqnarray}\label{eq:2d-bnd}
\int_0^T t\,\|u\|_{2}^2\,\ud t\le C
\end{eqnarray}
where the constant $C$ depends on $\Omega$, $\nu$, $\|u_0\|$ and $\int_0^T\|f\|^2\,\ud t$.
\end{lemma}


\begin{proof}

Let $A=\Pi(-\Delta)$ be the Stokes operator, where $\Pi$ is the orthogonal projection from
$L^2(\Omega)$ onto 
$$H=\{u\in L^2(\Omega):\nabla\cdot u=0\,\mbox{and } u\cdot{\bf n}=0\,\,\mbox{in the trace sense}\}$$
with ${\bf{n}}$ the outward normal vector on the boundary.
We take the inner product of (\ref{eq:nse}) with $t\,Au$.
Since $(\nabla p,Au)=0$, we can write
\begin{eqnarray*}
\fl\quad\frac{1}{2}\frac{\ud}{\ud t}(t\,\|A^{1/2}u\|^2)-\|A^{1/2}u\|^2&+\nu\, t\,\|Au\|^2 \\
&\le t\,\left((u\cdot\nabla)u,Au\right)+t\,(f,Au)\\
&\le c\,t\,\|u\|^{1/2}\,\|u\|_1\,\|u\|_2^{3/2}+t\,\|f\|\|Au\|\\
&\le c\,t\,\,\|u\|^2\,\|u\|_1^4+c\,t\,\|f\|^2+\frac{\nu}{2}t\,\|Au\|^2\\
&\le c\,t\,\,\|u\|^2\,\|u\|_1^2\,\|A^{1/2}u\|^2+c\,t\,\|f\|^2+\frac{\nu}{2}t\,\|Au\|^2,
\end{eqnarray*}
since $\|u\|_1\le c\,\|A^{1/2} u\|$ (see, for example, Robinson 2001, Proposition 6.18). Therefore
$$
\frac{\ud}{\ud t}(t\,\|A^{1/2}u\|^2)+\nu\,t\,\|Au\|^2 \le (c\,\|u\|^2\,\|u\|_1^2)\,t\,\|A^{1/2}u\|^2
+c\,t\,\|f\|^2+c\,\|u\|_1^2.
$$
Multiplying both sides of the above inequality with 
$E_2(t)=\exp\left(-c\,\int_0^t\|u\|^2\,\|u\|_1^2\,\ud s\right)$ we obtain
$$
\frac{\ud}{\ud t}(t\,E_2(t)\,\|A^{1/2}u\|^2)+\nu\,t\,E_2(t)\,\|Au\|^2\le c\,t\,\|f\|^2+c\,\|u\|_1^2.
$$
Now integrating the above inequality from $0$ to $T$ (which can be made rigorous using the Galerkin approximation) we obtain
$$
\int_0^T t\,\|u\|_2^2\,\ud t\,\le c\,\e^{c\,\|u\|_0^2}\,\int_0^T\|u\|_1^2+\|f\|^2\,\ud\tau,
$$
since $\|u\|_2\le c\,\|Au\|$. The result follows.
\end{proof}

This, by Theorem {\ref{thm:g-uniqueness}}, proves the uniqueness of the 
particle trajectories.


%
%
\subsection{The three-dimensional case}\label{sec:3d}

In this case when obtaining the bound on $\int_0^T \|u\|_{5/2}\,\ud\tau$
we will need to use the equivalence of 
$\|u\|_r$ and $\|\cA^{r/2}u\|$ where $\cA=-\Delta$ with Dirichlet boundary conditions, 
and $r$ is any non-negative real number. 
For $0\le r\le 2$, this equivalence follows from Fujiwara's (1967) characterization of $D(\cA^r)$.  
For $r>2$, to our knowledge, such a characterization does not exist.
However, noting that $\|u\|_r$ and $\|\cA^{r/2}u\|$ are equivalent when $r$ is a non-negative integer (Gilbarg \& Trudinger 1983, Robinson 2001), 
the result for positive real $r$ follows from an interpolation theorem of Lions \& Magenes (1972). We could not find this equivalence result stated explicitly 
in the literature and therefore we think it would be worthwhile to show it here.

We first note that we consider the fractional Sobolev spaces $H^s$ characterised by the following norm, for $s=m+\sigma$ with $m\in\mathbb{Z}$ and $0<\sigma<1$ (Adams 1975)
$$
\|u\|_s^2=\|u\|_m^2
+\sum_{|\alpha|=m}\int_{\Omega}\int_{\Omega}
    \frac{|D^\alpha u(x)-D^\alpha u(y)|^2}{|x-y|^{d+2\sigma}}\,\ud x\,\ud y.
$$
However, we assume $\Omega$ bounded with $\partial\Omega$ uniformly $C^1$-regular and therefore the above characterisation is equivalent to the definition based on the real interpolation method (Adams 1975, Theorem 7.48; Lions and Magenes 1972, Theorem 9.1 of Chapter 1).

\begin{lemma}\label{l:realp}
Let $\Omega$ be an open subset of $\R^d$, $d\ge 2$, with a sufficiently smooth boundary. Then for any non-negative real number $r$
\begin{equation}\label{eq:realp}
c_2\,\|\cA^{r/2}u\|\le\|u\|_r \le c_1\,\|\cA^{r/2}u\|, 
\quad \mbox{for all }\; u\in D(\cA^{r/2}).
\end{equation}
\end{lemma}
\begin{proof}

To prove $\|u\|_r \le c_1\,\|\cA^{r/2}u\|$, let $\{X_1,X_2\}$ and $\{Y_1,Y_2\}$ two pairs of normed spaces with $X_1$ and $Y_1$ dense subsets of 
$X_2$ and $Y_2$ respectively, with continuous injections. 
Define the interpolation 
space $[X_1,X_2]_{\theta}$ with $0<\theta\le 1$ as
$$
[X_1,X_2]_{\theta}=D(\Lambda^{1-\theta})
$$
where $\Lambda$ is a self adjoint, positive and unbounded operator in $X_2$, with 
domain $X_1$ and satisfying $(u,v)_{X_1}=(\Lambda u, \Lambda v)_{X_2}$, for any
$u,v\in X_1$. Also let
$$
[Y_1,Y_2]_{\theta}=D(S^{1-\theta})
$$
with $S$ having similar properties as $\Lambda$. The interpolation result of 
Lions \& Magenes (1972, Chapter I, Theorem 5.1) states that if a map $\pi$ is a continuous linear operator of $X_1$ into $Y_1$ and also of $X_2$
into $Y_2$, then it is a continuous linear operator from $[X_1,X_2]_{\theta}$ 
into $[Y_1,Y_2]_{\theta}$. 

For our purpose here, let $X_1=D(\cA^{m})$, $X_2=D(\cA^0)=H^0$, $Y_1=H^{2m}$ and $Y_2=H^0$, with $m$ 
a non-negative integer number.
Then $X_1$ and $Y_1$ are dense subsets of $X_2$ and $Y_2$ respectively with continuous
injections. 
We note that
for any real $r$ by definition of $H^r(\Omega)$ (Lions \& Magenes 1972),
we have 
$$
H^r(\Omega)=[H^{2m}(\Omega),H^0(\Omega)]_{\theta}\quad\mbox{with } \theta=1-r/(2m).
$$
Also for $X_1$ and $X_2$, letting $\Lambda=\cA^{m}$ implies that 
$$
D(\cA^{r/2})=[D(\cA^{m}),H^0(\Omega)]_{\theta}\quad\mbox{with } \theta=1-r/(2m).
$$
Now since for an integer $m$, $D(\cA^{m})\subset H^{2m}$, the
identity operator is a linear continuous operator from $D(\cA^{m})$ into $H^{2m}$
and also obviously from $H^0(\Omega)$ into $H^0(\Omega)$. Therefore by 
the result mentioned above it is a continuous operator of 
$D(\cA^{r/2})=[D(\cA^{m}),H^0(\Omega)]_{\theta}$  into
$H^r(\Omega)=[H^{2m}(\Omega),H^0(\Omega)]_{\theta}$, implying that
$D(\cA^{r/2})\subset H^r(\Omega)$ and therefore
$$
\|u\|_r\le c_1\,\|\cA^{r/2}u\|.
$$

It remains to show that $\|\cA^{r/2}u\|\le \|u\|_r$.
For any $r\ge 0$, there exists an integer $m\ge 0$ such that $r=2m+\hat{r}$ with real $0\le \hat{r}<2$.
Noting that by the result of Fujiwara (1967)
$$
\|A^{\hat{r}/2}u\|\,\le\, \|u\|_{\hat{r}},
$$
we can write
\begin{eqnarray*}
\|\cA^{r/2}u\|=\|\cA^{m+\hat{r}/2}u\|
&=\,\|\cA^{\hat{r}/2}\cA^{m}u\|\\
&\le\,\sum_{|\alpha_j=2|,\,1\le j\le m} \|D^{\alpha_1+\dots+\alpha_m}u\|_{\hat{r}}\\
&\le\,\|u\|_{2m+\hat{r}}=\|u\|_{r},
\end{eqnarray*}
and the result follows.

\end{proof}
%
%

Having (\ref{eq:realp}), we can also show that $u$, the local unique solution of (\ref{eq:nse}) with $u_0\in H^{1/2}$ 
and $f\in L^2(0,T;H^{1/2})$, is bounded in $C((0,T);D(\cA^{3/4}))$. Consider $\{w_1,w_2,\dots,w_m\}$ to be the first $m$ 
eigenfunctions of $\cA$ with corresponding eigenvalues $\{\lambda_1,\lambda_2,\dots,\lambda_m\}$. Then
since the Galerkin approximation $u_m\in \mathrm{span}\{w_1,\dots,w_m\}$ satisfies 
$\sum_{j=1}^m|\lambda_j|^{2r}|(u_m,w_j)|^2<\infty$ 
for any finite $r$ and therefore is in $D(\cA^r)$,
it can be shown in a similar way to the proof of Lemma 4.2 of Temam (1995) that 
$u\in C((0,T);D(\cA^{3/4}))$ and therefore $\|\cA^{5/4}u\|<\infty$ for almost every
$t\in (0,T)$. The fact that $u(t)\in D(\cA^{5/4})$ for almost every $t>0$ is used in the proof of
the next lemma where we not only require enough Sobolev regularity for $u$, 
but also need it to be in $D(\cA^{5/4})$ to be able to use (\ref{eq:realp}).\\

We can now show the bound on $\int_0^T t\,\|u\|_{5/2}^2\,\ud t$:\\


\begin{lemma}\label{lem:3d}
Let $\Omega$ be an open bounded subset of $\R^3$ with a sufficiently smooth boundary.
Consider $u\in L^\infty(0,T;H^{1/2}(\Omega))\cap L^2(0,T;H^{3/2}(\Omega))$ to be a unique solution of 
the three-dimensional Navier-Stokes equations (\ref{eq:nse}) with
$u_0\in H^{1/2}(\Omega)$ and $f\in L^2(0,T;H^{1/2}(\Omega))$. 
Then 
\begin{eqnarray}\label{eq:3d-bnd}
\int_0^T t\,\|u\|_{5/2}^2\,\ud t \le C
\end{eqnarray}
where $C$ depends on $\nu$, $\Omega$, $\|u_0\|_{1/2}$ and $\int_0^T\|f\|_{1/2}^2\,\ud t$.
\end{lemma}


\begin{proof}

We take the inner product of (\ref{eq:nse}) with $t\,\cA^{3/2} u$ to obtain
\begin{eqnarray}\label{eq:energy3d}
\fl\quad\frac{1}{2}\frac{\ud}{\ud t}(t\,\|\cA^{3/4}u\|^2)
-\|\cA^{3/4}u\|^2+\nu\,t\,\|\cA^{5/4}u\|^2
&\le t\,|\left((u\cdot \nabla)u ,\cA^{3/2}u\right)|\nonumber\\
&\;+t\,|(\nabla p,\cA^{3/2}u)|+t\,|(f,\cA^{3/2}u)|.
\end{eqnarray}
The highest derivative exponent in the right-hand side is bigger than that of the left-hand side. Therefore we 
need to integrate by parts in the right-hand side, which is why in the three-dimensional case we have to take
the inner product of (\ref{eq:nse}) with $\cA^{3/2} u$ rather than $A^{3/2} u$ 
(noting that $A=\Pi(-\Delta)\neq -\Delta$ in the bounded domains).

For the second term in the right-hand side of (\ref{eq:energy3d}), by integration by parts and appropriate use of 
the Sobolev embeddings, we obtain (noting that $u|_{\partial\Omega}=0$)
\begin{eqnarray*}
\fl|\left((u\cdot \nabla)u ,\cA^{3/2}u\right)|\\
\fl\qquad\le \|\,|Du|\,|Du|\,|D\cA^{1/2}u|\,\|_{L^1(\Omega)}
          +\|\,|u|\,|D^2u|\,|D\cA^{1/2}u|\,\|_{L^1(\Omega)}\\
\fl\qquad\le \|Du\|_{L^3(\Omega)}^2\,\|D\cA^{1/2}u\|_{L^3(\Omega)}
          +\|u\|_{L^6(\Omega)}\,\|D^2u\|\,\|D\cA^{1/2}u\|_{L^3(\Omega)}\\
\fl\qquad\le c\,\|u\|_{{3/2}}^2\,\|u\|_{{5/2}}+c\,\|u\|_1\,\|u\|_2\,\|u\|_{5/2}\\
\fl\qquad\le c\,\|u\|_{{3/2}}^2\,\|u\|_{{5/2}}
                 +c\,\|u\|_{1/2}^{1/2}\,\|u\|_{3/2}\,\|u\|_{5/2}^{3/2}.
\end{eqnarray*}
For the term containing the pressure, again by integration by parts, we can write
$$
|(\nabla p,\cA^{3/2}u)|
\le \|\,|D^2p|\,|D\cA^{1/2}u|\,\|_{L^1(\Omega)}+\|\,|Dp|\,|D\cA^{1/2}u|\,\|_{L^1(\partial\Omega)}.
$$
Since applying the divergence operator to (\ref{eq:nse}) gives (assuming that $f$ is also divergence-free)
$$
-\Delta p=\nabla\cdot\left( (u\cdot\nabla)u \right)=\sum_{i,j=1}^3\partial_ju_i\partial_iu_j
$$
and also $\|\cA p\|_{L^q}\le\|D^2 p\|_{L^q}\le c\,\|\cA p\|_{L^q}$ for any $q>1$ 
(Gilbarg \& Trudinger 1983, Lemma 9.17), we can write
\begin{eqnarray*}
\|\,|D^2p|\,|D\cA^{1/2}u|\,\|_{L^1(\Omega)}
&\le\, \|D^2p\|_{L^{3/2}(\Omega)}\,\|D\cA^{1/2}u\|_{L^3(\Omega)}\\
&\le\, c\,\|Du\|_{L^3(\Omega)}^2\,\|u\|_{5/2}\,\le\, c\,\|u\|_{3/2}^2\,\|u\|_{5/2}
\end{eqnarray*}
and
\begin{eqnarray*}
\|\,|Dp|\,|D\cA^{1/2}u|\,\|_{L^1(\partial\Omega)}
&\le\, c\,\|p\|_{H^1(\partial\Omega)}\,\|u\|_{H^2(\partial\Omega)}
           \,\le\,c\,\|p\|_{3/2}\,\|u\|_{5/2}\\
&\le\, c\,\|D^2p\|_{L^{3/2}(\Omega)}\,\|u\|_{5/2}
           \,\le\, c\,\|Du\|_{L^3(\Omega)}^2\,\|u\|_{5/2}\\
&\le\, c\,\|u\|_{3/2}^2\,\|u\|_{5/2}.
\end{eqnarray*}
Substituting these in (\ref{eq:energy3d}) and using (\ref{eq:realp}) (which holds for $r\le 5/2$ and almost every $t>0$, 
since $u(t)\in D(\cA^{5/4})$ for almost every $t>0$) we conclude that
\begin{eqnarray*}
\fl\qquad\frac{\ud}{\ud t}(t\,\|\cA^{3/4}u\|^2)+\nu\,t\,\|\cA^{5/4}u\|^2\\
\fl\qquad\qquad\le c\,\|u\|_{3/2}^2\,(1+\|u\|_{1/2}^2)\,t\,\|\cA^{3/4} u\|^2
+c\,\|u\|_{3/2}^2+c\,t\,\|\cA^{1/4} f\|^2.
\end{eqnarray*}
Let
$E_3(t)=\exp\left(-c\int_0^t \|u\|_{3/2}^2\,(1+\|u\|_{1/2}^2)\,\ud s  \right)$ and multiply both sides of the above inequality by $E_3(t)$ to obtain
$$
\frac{\ud}{\ud t}(t\,E_3(t)\,\|\cA^{3/4}u\|^2)+\nu\,t\,E_3(t)\,\|\cA^{5/4}u\|^2
\le c\,\|u\|_{3/2}^2+c\,t\,\|\cA^{1/4} f\|^2.
$$
Integrating the above inequality between $0$ and $t$ (noting that this can be made rigorous using the Galerkin approximations of $u$) gives 
$$
\int_0^T t\,\|\cA^{3/4}u\|^2\,\ud t\,
\le \, \frac{c}{E_3(T)}\, \int_0^T(\|u\|_{3/2}^2+\|f\|_{1/2}^2)\,\ud t,
$$
and the result follows. 

\end{proof}

We therefore conclude the uniqueness of particle trajectories in the three-dimensional case as well.


\section{Conclusion}

We considered the two- and three-dimensional Navier-Stokes equations with  the initial conditions that have minimal Sobolev regularity required to give rise to a unique solution,
and presented a much simpler proof than that of
Chemin \& Lerner (1995) for the uniqueness of particle trajectories associated to such solutions.

For $u$ a particular weak solution 
of the three-dimensional Navier-Stokes equations, it is shown by Foias, Guillop\'e \& Temam (1985) that at least one continuous solution of (\ref{eq:ode}) exists.
The uniqueness of these solutions however is not known. 
Finding an extra condition on a weak solution $u$ that can lead to the uniqueness of
the solution of (\ref{eq:ode}) is the subject of the recent work 
by Robinson \& Sadowski (2008).

We note that the above result and also the continuity and differentiability properties of the
fluid particle trajectories with respect to the initial velocity field, are useful in showing that the posterior measure in certain data assimilation problems is well-defined 
(Cotter {\it et al}, 2008).  It can be shown (Cotter {\it et al}, 2008) that the two-dimensional trajectories are Lipschitz continuous and also differentiable with respect to an initial velocity field which is $H^s$-regular for some $s>0$. 
The continuity of the trajectories with respect to an only $H^{(d/2)-1}$-regular initial condition follows from a general continuity result for ordinary differential equations proved in Dashti \& Robinson (2008).

\ack
We would like to thank Andrew Stuart for suggesting the problem of uniqueness of particle trajectories of weak solutions of the two-dimensional equations, and Isabelle Gallagher for helpful discussion. M.D. was partially supported by a Warwick Postgraduate Research Fellowship, and EPSRC grant ER/F050798/1.
J.C.R. was partially supported by a Royal Society University Research Fellowship and, 
an EPSRC Leadership Fellowship EP/G007470/1. We would also like to thank the referees for their helpful comments.

\section*{References}
\begin{harvard}

\item[] Adams, R.A. 1975, {\it Sobolev spaces.} Pure and Applied Mathematics, Vol. 65. Academic Press, New York-London, 1975.

\item[] Chemin J.-Y. \& Lerner N. 1995, Flot de champs de vecteurs non lipschitziens et \'equations de Navier-Stokes. 
{\it J. Differential Equations}, {\bf 121}, 314--328 

\item[] Cotter S., Dashti M., Robinson J.C. and Stuart A.M. 2008, Data assimilation problems in fluid mechanics: Bayesian formulation in function space. {\it Submitted}.

\item[] Dashti M. \& Robinson J.C. 2008, The uniqueness of Lagrangian trajectories in Navier-Stokes flows, to appear in Robinson J.C. and Rodrigo J.L. {\it Partial Differential Equations and Fluid Mechanics} (Warwick, May 2007), to be published by Cambridge University Press, Cambridge.

\item[] Foias C., Guillop\'e C. \& Temam R. 1985, Lagrangian representation of a flow. \textit{Journal of Differential Equations} \textbf{57}, 440--449.

\item[] Fujita H. \& Kato T. 1964, On the Navier-Stokes initial value problem I. {\it Arch. Rat. Mech. Anal.}, {\bf 16}, 269--315.

\item[] Fujiwara D. 1967, Concrete characterization of the domains of fractional powers of some elliptic differential operators of the second order. {\it Proc. Japan Acad.}, {\bf 43}, 82--86.

\item[] Gilbarg D. \& Trudinger N.S. 1983, {\it Elliptic partial differential equations of second order}. Springer-Verlag, Berlin.

\item[] Ladyzhenskaya O.A. 1958, Solution in the large to the boundary-value problem for the Navier-Stokes equations in two-space variables. {\it Soviet Physics. Dokl.} 123, 1128-1131.

\item[] Leray J. 1933, Etude de diverses \'equations int\'egrates non lineaires et de quelques probl\'emes que pose l'hydrodynamique. {\it J. Math. Pures Appl.}, 12, 1-82.

\item[] Lions J.-L. \& Magenes E. 1972, {\it Non-homogeneous boundary value problems and applications.} Vol. I, Springer-Verlag, New York.

\item[] Lions, J-L. \& Prodi, G. 1959, Un thŽorme d'existence et unicitŽ dans les Žquations de Navier-Stokes en dimension 2. (French) {\it C. R. Acad. Sci. Paris}  248, 3519--3521.

\item[] Priestley, H. A. 1997, {\it Introduction to integration.} Oxford Science Publications. The Clarendon Press, Oxford University Press, New York.

\item[] Robinson J.C. 2001, {\it Infinite-dimensional dynamical systems.} Cambridge Texts in Applied
Mathematics, Cambridge University Press, Cambridge.

\item[] Robinson J.C., Sadowski W. 2008, A criterion for uniqueness of Lagrangian trajectories for weak solutions of the 3d Navier-Stokes equations. {\it Submitted}.

\item[] Temam R. 1995, {\it Navier-Stokes Equations and Nonlinear Functional Analysis}. SIAM, Philadelphia.

\item[] Temam R. 1977, {\it Navier-Stokes Equations}. AMS Chelsea Publishing, Providence, RI.

\item[] Zuazua, E. 2002, Log-Lipschitz regularity and uniqueness of the flow for a field in $(W\sp {n/p+1,p}\sb {\rm loc}(\Bbb R\sp n))\sp n$. {\it C. R. Math. Acad. Sci. Paris}, {\bf 335}, no. 1, 17--22.

\end{harvard}

\end{document}